\documentclass{amsart}
\usepackage{amssymb}
\usepackage{a4} 
\begin{document}
\newtheorem{theorem}{Theorem}[section]
\newtheorem{lemma}[theorem]{Lemma}
\newtheorem{remark}[theorem]{Remark}
\newtheorem{definition}[theorem]{Definition}
\newtheorem{corollary}[theorem]{Corollary}
\newtheorem{example}[theorem]{Example}
\def\Pspan{\operatorname{Span}}
\def\Rank{\operatorname{Rank}}
\def\BB{\mathcal{B}}
\def\RR{\mathfrak{R}}
\def\qedbox{\hbox{$\rlap{$\sqcap$}\sqcup$}}
\makeatletter
  \renewcommand{\theequation}{%
   \thesection.\alph{equation}}
  \@addtoreset{equation}{section}
 \makeatother
\title[Curvature homogeneous manifolds]
{Curvature homogeneous pseudo-Riemannian manifolds which are not locally homogeneous}
\author{Corey Dunn and Peter B. Gilkey${}^1$}
\thanks{${}^1$ Research partially supported by the Max Planck Institute (Leipzig)}
\begin{address}{Mathematics Department, University of Oregon, Eugene Or 97403 USA}
\end{address}
\begin{email}{cdunn@darkwing.uoregon.edu and gilkey@darkwing.uoregon.edu}\end{email}
\begin{abstract} We construct a family of balanced signature pseudo-Riemannian
manifolds, which arise as hypersurfaces in flat space, that are curvature homogeneous,
that are modeled on a symmetric space, and that are not locally homogeneous.
\end{abstract}
\keywords{curvature homogeneous, balanced signature, hypersurfaces}
\subjclass[2000]{Primary 53C50}
\maketitle

\section{introduction}\label{sect1} Let $R$ be the Riemann curvature tensor of a
pseudo-Riemannian manifold $(M,g)$ of signature $(p,q)$. Following Kowalski, Tricerri, and
Vanhecke \cite{KTV91,KTV92}, we say that $(M,g)$ is {\it curvature homogeneous} if given any two points
$P,Q\in M$, there is a linear isomorphism $\Psi:T_PM\rightarrow T_QM$ such that $\Psi^*g_Q=g_P$
and such that $\Psi^*R_Q=R_P$; this notion has also been called $0$ curvature homogeneous when
considering a similar condition for the higher covariant derivatives of the curvature tensor.

Similarly, $(M,g)$ is said to be {\it locally homogeneous} if given any two points $P$ and $Q$,
there are neighborhoods $U_P$ and $U_Q$ of $P$ and $Q$, respectively, and an isometry
$\psi:U_P\rightarrow U_Q$ such that $\psi P=Q$. Taking $\Psi:=\psi_*$ shows that locally
homogeneous manifolds are curvature homogeneous. The somewhat surprising fact is that the converse
fails -- there are curvature homogeneous manifolds which are {\bf not} locally homogeneous. 

There is by now an extensive literature on the subject in the Riemannian setting, see, for example, the
discussion in \cite{BK96,BV98,K98,T97,T88,V91}. There are also a number of papers in the Lorentzian setting
\cite{BM00, BV97,CLPT90} and also in the affine setting \cite{KO99,O96}. There are, however, almost no papers in
the higher dimensional setting -- and those that exist appear in the study of $4$ dimensional neutral signature
Osserman manifolds, see, for example, \cite{BCR01,GKV02}. In this brief note, we exhibit a
family of examples in signature $(p,p)$ for any $p\ge3$ which are curvature homogeneous but not locally homogeneous;
this family first arose in the study of Szab\'o Osserman IP Pseudo-Riemannian manifolds \cite{GIZ02,GIZ03}.

Let $( x, y)=(x_1,...,x_p,y_1,...,y_p)$ be the usual coordinates on
$\mathbb{R}^{2p}$. Let $f(x)$ be a smooth function on an open subset $\mathcal{O}\subset\mathbb{R}^p$. We
define a non-degenerate pseudo-Riemannian metric $g_f$ of balanced signature $(p,p)$ on
$M:=\mathcal{O}\times\mathbb{R}^p$ by:
\begin{equation}\label{eqn-1.a}
g_f(\partial_i^x,\partial_j^x)=\partial_i^xf\cdot\partial_j^xf,\quad
g_f(\partial_i^x,\partial_i^y)=\delta_{ij},\quad\text{and}\quad
g_f(\partial_i^y,\partial_j^y)=0\,.
\end{equation}
This is closely related to the so called `deformed complete lift' of a metric on 
$\mathcal{O}$ to $T\mathcal{O}$, see, for example, the discussion in
\cite{BCGHV98,Ka00,Opr89}.

The pseudo-Riemannian manifold $(M,g_f)$ arises as a hypersurface in a flat space. Let
$\{\vec u_1,...,\vec u_p,\vec v_1,...,\vec v_p,\vec w_1\}$ be a basis for a vector space $W$. Define an inner product
$\langle\cdot,\cdot\rangle$ of signature
$(p,p+1)$ on $W$ by setting
$$\begin{array}{lll}
     \langle\vec u_i,\vec u_j\rangle=0,&
     \langle\vec u_i,\vec v_j\rangle=\delta_{ij},
    &\langle\vec v_i,\vec v_j\rangle=0,\\
     \langle\vec u_i,\vec w_1\rangle=0,&
     \langle\vec v_i,\vec w_1\rangle=0,&
     \langle\vec w_1,\vec w_1\rangle=1\,.\vphantom{\vrule height 11pt}
\end{array}$$
Let $F(x,y)=x_1\vec u_1+...+x_p\vec u_p+y_1\vec v_1+...+y_p\vec v_p+f(x)\vec w_1$ define an embedding of $M$ in $W$.
Then $g_f$ is the induced metric on the embedded hyper surface. The normal $\nu$ to the hypersurface is given by setting
$\nu=\vec w_1-\partial_1^xf\ \vec v_1-...-\partial_p^xf\ \vec v_p$. Thus the second fundamental form $L_f$ of the
embedding is given by the Hessian
$$
L_f(\partial_i^x,\partial_j^x)=\partial_i^x\partial_j^xf,\quad
L_f(\partial_i^x,\partial_j^y)=0,\quad\text{and}\quad
L_f(\partial_i^y,\partial_j^y)=0\,.
$$
We define distributions
$$
\mathcal{X}:=\Pspan\{\partial_1^x,...,\partial_p^x\}\quad\text{ and }\quad
  \mathcal{Y}:=\Pspan\{\partial_1^y,...,\partial_p^y\}.
$$
We then have $L(Z_1,Z_2)=0$ if $Z_1\in\mathcal{Y}$ or $Z_2\in\mathcal{Y}$ so the restriction $L_f^{\mathcal{X}}$
of  $L$ to the distribution $\mathcal{X}$ carries the essential information. The following is the main result of this
paper:
\begin{theorem}\label{thm-1.1} 
If the quadratic form $L_f^{\mathcal X}$ is positive definite, then $(M,g_f)$ is
curvature homogeneous. Furthermore, if $p\ge3$, then $(M,g_f)$ is not locally homogeneous for generic $f$.
\end{theorem}

As noted above, these manifolds first arose in an entirely different setting. 
Let $R$ be the Riemann curvature tensor of a pseudo-Riemannian manifold $(M,g)$. Let $\nabla R$ be the covariant
derivative of $R$. Let
$J$, $S$, and
$\mathcal{R}$ be the associated Jacobi operator, Szab\'o operator, and skew-symmetric curvature operator, respectively. 
Let $X\in TM$ and let $\{Y,Z\}$ be an oriented orthonormal basis for an oriented spacelike or timelike $2$ plane $\pi$. 
These operators are defined by the identities:
\begin{eqnarray*}
&&g(J(X)U,V)=R(U,X,X,V),\\
&&g(S(X)U,V)=\nabla R(U,X,X,V;X),\\
&& g(\mathcal{R}(\pi)U,V)=R(Y,Z,U,V)
\end{eqnarray*}
Stanilov and Videv \cite{SV98} have defined a {\it higher order Jacobi operator} by setting
$$J(\pi):=g(X_1,X_1)J(X_1)+...+g(X_\ell,X_\ell)J(X_\ell)$$
where $\{X_1,...,X_\ell\}$ is any orthonormal basis for a non-degenerate subspace $\pi\subset TM$. 

\begin{definition}\rm Let $(N,g)$ be a pseudo-Riemannian manifold. Then $(N,g)$ is\begin{enumerate}
\item {\it spacelike Jordan Osserman} (resp. {\it timelike Jordan Osserman}) if the Jordan
normal form of
$J(X)$ is constant on the bundle of unit spacelike (resp. unit timelike) vectors.
\item {\it spacelike Szab\'o} (resp. {\it timelike Szab\'o}) if the eigenvalues of $S(X)$ are constant on the bundle of
unit spacelike (resp. unit timelike) vectors.
\item {\it spacelike Jordan IP} (resp. {\it timelike Jordan IP}) if the Jordan normal form of $\mathcal{R}(\pi)$
is constant on the Grassmannian of oriented spacelike (resp. timelike) $2$ planes in $TM$.
\item {\it Jordan Osserman} of type $(r,s)$ if the Jordan normal form of $J(\pi)$ is constant on the Grassmannian of
non-degenerate subspaces of type $(r,s)$ in $TM$.
\end{enumerate} 
\end{definition}

The spectral geometry of the
Jacobi operator, of the skew-symmetric curvature operator, and of the Szab\'o operator were first considered in the
Riemannian setting by Osserman
\cite{Os97}, by Ivanova and Stanilov \cite{IS95}, and by Szab\'o \cite{Sz91}, respectively. We refer to \cite{G02} for
further details. The manifolds
$(M,g_f)$ provide examples of these manifolds. We refer to \cite{GIZ02,GIZ03} for the proof of the following result:

\begin{theorem}\label{gronk1.2} If the quadratic form $L_f^{\mathcal X}$ is positive definite, then $(M,g_f)$ is
spacelike Jordan Osserman, timelike Jordan Osserman, spacelike Szab\'o, timelike Szab\'o, spacelike Jordan IP,
and timelike Jordan IP. Furthermore $(M,g_f)$ is Jordan Osserman of types $(r,0)$, $(0,r)$, $(p-r,p)$, and
$(p,p-r)$,and is not Jordan Osserman of type $(r,s)$ otherwise.
\end{theorem}

We note there are no known Jordan Szab\'o manifolds which are not symmetric.

Here is a brief guide to the paper. In Section \ref{Sect2}, we determine the tensors $R_f$ and $\nabla R_f$ which are
defined by the metric
$g_f$ and show $(M,g_f)$ is curvature homogeneous. In Section \ref{Sect3}, we complete the proof of Theorem
\ref{thm-1.1} by showing that $(M,g_f)$ is not locally homogeneous for generic $f$. We conclude in Remark
\ref{gronk3.3} by showing the `model space' for the curvature tensor for $(M,g_f)$ is that of a symmetric space.

It is a pleasant task to acknowledge helpful conversations on this subject with Prof. Garc\'{\i}a--R\'{\i}o and
with Prof. J.H. Park. The paper is dedicated to Professor Vanhecke. The second author has had the priviledge of knowing
Professor Vanhecke for a number of years and owes Professor Vanhecke a profound debt of gratitude not only for many
useful mathematical discussions but also for wise counsel on a number of subjects.

\section{The tensors $R_f$ and $\nabla R_f$}\label{Sect2} We begin the proof of Theorem \ref{thm-1.1} by determining
$R_f$ and $\nabla R_f$.

\begin{lemma}\label{gronk2.1}
 Let $Z_1,...$ be coordinate vector fields on $M:=\mathcal{O}\times\mathbb{R}^p$. Let the metric $g_f$ be given by
Equation {\rm(\ref{eqn-1.a})}. Then:
\begin{enumerate}
\item $\nabla_{Z_1}Z_2=0$ if $Z_1\in\mathcal{Y}$ or if $Z_2\in\mathcal{Y}$;
\item $R(Z_1,Z_2,Z_3,Z_4)=L(Z_1,Z_4)L(Z_2,Z_3)-L(Z_1,Z_3)L(Z_2,Z_4)$. This vanishes if one of the $Z_i\in\mathcal{Y}$
for $1\le i\le4$;
\item $\nabla R(Z_1,Z_2,Z_3,Z_4;Z_5)=Z_5\{R(Z_1,Z_2,Z_3,Z_4)\}$. This vanishes if one of the $Z_i\in\mathcal{Y}$ for
$1\le i\le5$.
\end{enumerate}\end{lemma}

\begin{proof} We have
$$(\nabla_{Z_1}Z_2,Z_3)=\textstyle{\frac12}\{Z_2g_f(Z_1,Z_3)+Z_1g_f(Z_2,Z_3)-Z_3g_f(Z_1,Z_2)\}\,.$$
This vanishes if any of the $Z_i\in\mathcal{Y}$. Assertion (1) now follows. Let $g_{ij}^x:=g(\partial_i^x,\partial_j^x)$
and let $\Gamma_{ijk}^x:=\frac12(\partial_i^xg_{jk}^x+\partial_j^xg_{ik}^x-\partial_k^xg_{ij}^x)$. We adopt the Einstein
convention and sum over repeated indices to see
$$\nabla_{\partial_i^x}\partial_j^x=\Gamma_{ijk}^x\partial_k^y,\quad
\nabla_{\partial_i^x}\partial_j^y=\nabla_{\partial_j^y}\partial_i^x=0,\quad\text{and}\quad
\nabla_{\partial_i^y}\partial_j^y=0\,.
$$
It now follows that $R(Z_1,Z_2,Z_3,Z_4)=0$ if any of the $Z_i\in\mathcal{Y}$. Furthermore
$$R(\partial_i^x,\partial_j^x,\partial_k^x,\partial_l^x)=\partial_i\Gamma_{jkl}^x-\partial_j\Gamma_{ikl}^x\,.$$
Assertion (2) now follows; this also, of course, follows from the classical formula which expresses the curvature
tensor of a hypersurface in flat space in terms of the second fundamental form.

Since $\nabla_{Z_5}Z_i\in\mathcal{Y}$ and since $R(\cdot,\cdot,\cdot,\cdot)$ vanishes if any of the entries belong to
$\mathcal{Y}$, Assertion (3) follows from Assertion (2).
\end{proof}

\medbreak We show that $(M,g_f)$ is curvature homogeneous by showing the following result:
\begin{lemma}\label{gronk2.2} Let $P\in M$. Assume $L_f^{\mathcal{X}}$ is positive definite. Then there exists a basis
$\{X_1,...,X_p,Y_1,...,Y_p\}$ for $T_PM$ so that:
\begin{enumerate}
\item $g_f(X_i,X_j)=0$, $g_f(X_i,Y_j)=\delta_{ij}$, and $g_f(Y_i,Y_j)=0$.
\item $R_f(X_i,X_j,X_k,X_l)=\delta_{il}\delta_{jk}-\delta_{ik}\delta_{jl}$.
\item $R_f(\cdot,\cdot,\cdot,\cdot)=0$ if any
of the entries is one of the vector fields $\{Y_1,...,Y_p\}$.
\end{enumerate}
\end{lemma}

\begin{proof} Fix $P\in M$. We diagonalize the quadratic form $L_f^{\mathcal{X}}$ at $P$ to choose tangent vectors $\bar
X_i=a_{ij}\partial_j^x\in T_PM$ so that
$L(\bar X_i,\bar X_j)=\delta_{ij}$. Let $\bar Y_i:=a^{ji}\partial_j^y$ where $a^{ij}$ is the inverse matrix. Then
\begin{eqnarray*}
&&g_f(\bar X_i,\bar Y_j)=a_{ik}a^{\ell j}g_f(\partial_k^x,\partial_\ell^y)=a_{ik}a^{kj}=\delta_{ij},\\
&&g_f(\bar Y_i,\bar Y_j)=0,\\
&&R_f(\bar X_i,\bar X_j,\bar X_k,\bar X_\ell)=\delta_{i\ell}\delta_{jk}-\delta_{ik}\delta_{j\ell},
\end{eqnarray*}
and $R_f(\cdot,\cdot,\cdot,\cdot)=0$ if any entry is $\bar Y_i$. We define
$$X_i:=\bar X_i-\textstyle\frac12g_f(\bar X_i,\bar X_j)\bar Y_i\quad\text{and}\quad Y_i:=\bar Y_i$$
to ensure $g_f(X_i,X_j)=0$. It is immediate that the frame $\{X_1,...,X_p,Y_1,...,Y_p\}$ satisfies the normalizations of
the Lemma.
\end{proof}

\section{Homogeneity}\label{Sect3} We begin our discussion with a technical observation. Let
$V$ be a finite dimensional real vector space. A
$4$ tensor $R\in\otimes^4V^*$ is said to be an {\it algebraic curvature tensor} if it satisfies the symmetries of the
Riemann curvature tensor, i.e. if:
\begin{eqnarray*}
&&R(\vec v_1,\vec v_2,\vec v_3,\vec v_4)=-R(\vec v_2,\vec v_1,\vec v_3,\vec v_4)
=R(\vec v_3,\vec v_4,\vec v_1,\vec v_2)\quad\text{and}\\
&&R(\vec v_1,\vec v_2,\vec v_3,\vec v_4)+R(\vec v_2,\vec v_3,\vec v_1,\vec v_4)
 +R(\vec v_3,\vec v_1,\vec v_2,\vec v_4)=0\,.
\end{eqnarray*}
If $\phi$ is a symmetric bilinear form on $V$, then we may
define an algebraic curvature tensor $R_\phi$ on $V$ by setting:
$$R_\phi(\vec v_1,\vec v_2,\vec v_3,\vec v_4):=
  \phi(\vec v_1,\vec v_4)\phi(\vec v_2,\vec v_3)-\phi(\vec v_1,\vec v_3)\phi(\vec v_2,\vec v_4)\,.$$

\begin{lemma}\label{gronk3.1} Let $\phi_1$ and $\phi_2$ be symmetric positive definite bilinear forms on a vector
space $V$ of dimension at least $3$. If
$R_{\phi_1}=R_{\phi_2}$, then $\phi_1=\phi_2$.
\end{lemma}
We note that Lemma \ref{gronk3.1} fails if $\dim V\le 2$.

\begin{proof} Since $\phi_1$ is positive definite, we can diagonalize $\phi_2$ with respect to $\phi_1$ and choose a
basis $\{\vec e_1,...,\vec e_r\}$ for $V$ so that $\phi_1(\vec e_i,\vec e_j)=\delta_{ij}$ and so that
$\phi_2(\vec e_i,\vec e_j)=\lambda_i\delta_{ij}$. If $i\ne j$, then
\begin{eqnarray}\nonumber
1&=&\phi_1(\vec e_i,\vec e_i)\phi_1(\vec e_j,\vec e_j)-\phi_1(\vec e_i,\vec e_j)\phi_1(\vec e_i,\vec e_j)=R_{\phi_1}(\vec e_i,\vec e_j,\vec e_j,\vec e_i)
   \label{oink3.a}\\
&=&R_{\phi_2}(\vec e_i,\vec e_j,\vec e_j,\vec e_i)=\phi_2(\vec e_i,\vec e_i)\phi_2(\vec e_j,\vec e_j)-\phi_2(\vec e_i,\vec e_j)\phi_2(\vec e_i,\vec e_j)\\
&=&\lambda_i\lambda_j\,.\nonumber
\end{eqnarray}
Since $r\ge3$, we can choose $k$ so $\{i,j,k\}$ are distinct indices. By Equation (\ref{oink3.a}),
$1=\lambda_i\lambda_k=\lambda_j\lambda_k$ so $\lambda_i=\lambda_j$ for all $i,j$. Since
$1=\lambda_i\lambda_j=\lambda_i^2$ and since
$\phi_2$ is positive definite, $\lambda_i=1$ for all $i$ and hence $\phi_1=\phi_2$.
\end{proof}

We say that $\BB:=(X_1,...,X_p,Y_1,...,Y_p)$
is an {\it admissible} basis for $T_PM$ if $\BB$ satisfies the normalizations of Lemma \ref{gronk2.2}. We can now
define a useful invariant:

\begin{lemma}\label{gronk3.2} Suppose that $L_f^{\mathcal{X}}$ is positive definite. Let $P\in M$ and let
$\BB$ be an admissible basis for $T_PM$. Let
$\alpha_f(P,\mathcal{B}):=\sum_{i,j,k,l,n}\nabla R_f(X_i,X_j,X_k,X_l;X_n)(P)^2$.
\begin{enumerate}
\item $\alpha_f(P,\mathcal{B})$ is independent of the particular admissible basis $\BB$ which is chosen.
\item If $(M,g_f)$ is locally homogeneous, then $\alpha_f$ is the constant function.
\end{enumerate}
\end{lemma}

\begin{proof} The distribution
$\mathcal{Y}$ is invariantly defined being characterized by the property:
$$\mathcal{Y}_P=\{Y\in T_PM:R(Z_1,Z_2,Z_3,Y)=0\quad\text{for all}\quad Z_i\in T_PM\}\,.$$
The subspace $\mathcal{X}$ on the other hand is not invariantly defined. Denote the standard projection by $\pi$ from
$T_PM$ to $T_PM/\mathcal{Y}_P$. As 
$$L(\cdot,\cdot)=0,\quad
R_f(\cdot,\cdot,\cdot,\cdot)=0\quad\text{and}\quad\nabla R_f(\cdot,\cdot,\cdot,\cdot;\cdot)=0$$
if any entry belongs to
$Y$, these tensors induce corresponding structures $\bar L_f$, $\bar R_f$, and $\RR_f$ on $T_PM/\mathcal{Y}_P$ so that
$$L_f=\pi^*\bar L_f,\quad R_f=\pi^*\bar R_f,\quad\text{and}\quad \nabla R_f=\pi^*\RR_f\,.$$
If $\mathcal{B}$ is an admissible basis, then we may define a quadratic form $\phi_{\mathcal{B}}$ on
$T_PM/\mathcal{Y}_P$ by requiring that
$\{\pi X_1,...,\pi X_p\}$ is orthonormal with respect to this quadratic form. We then have $\bar
R_f=R_{\phi_{\mathcal{B}}}$. By Lemma \ref{gronk3.1}, $\phi=\phi_{\mathcal{B}}$ is independent of the particular
basis chosen and is invariantly defined. This defines a positive definite inner product on $T_PM/\mathcal{Y}_P$ which
we use to raise and lower indices and to contract tensors. The invariant $\alpha$ is then given by
$||\RR_f||_\phi^2$ and is invariantly defined. Since the structures involved are preserved by isometries, the Lemma
now follows. What we have done, of course, is to prove that the second fundamental form is preserved by a local
isometry of $(M,g_f)$ in this setting.
\end{proof}

\medbreak\noindent{\bf Proof of Theorem \ref{thm-1.1}}. In light of Lemma \ref{gronk3.2}, to complete the proof of
Theorem
\ref{thm-1.1}, it suffices to construct
$f$ so that
$\alpha_f$ is constant on no open subset of $\mathbb{R}^p$; the fact that such $f$ are generic will then follow using
standard arguments. Let $f_{;i}=\partial_i^xf$, $f_{;ij}:=\partial_i^x\partial_j^xf$, and so forth. We use Lemma
\ref{gronk2.1} to see:
\begin{eqnarray*}
&&R(\partial_i^x,\partial_j^x,\partial_k^x,\partial_l^x)=f_{;il}f_{;jk}-f_{;ik}f_{;jl},\\
&&\nabla R(\partial_i^x,\partial_j^x,\partial_k^x,\partial_l^x;\partial_n^x)
   =\partial_n^x\{f_{;il}f_{;jk}-f_{;ik}f_{;jl}\}\,.
\end{eqnarray*}
Let $\Theta=\Theta(x_1)$ be a smooth function on $\mathbb{R}$ so that $|\Theta_{;11}|<1$. Set
$$f(x):=\textstyle\frac12\{x_1^2+...+x_p^2\}+\Theta(x_1)\,.$$
 We may then compute, up to the usual $\mathbb{Z}_2$
symmetries, that the non-zero components of $R_f$ and of $\nabla R_f$ are:
$$\begin{array}{ll}
R_f(\partial_1^x,\partial_i^x,\partial_i^x,\partial_1^x)=1+\Theta_{;11}&\text{for}\quad 2\le i\le p,\\
R_f(\partial_i^x,\partial_j^x,\partial_j^x,\partial_i^x)=1&\text{for}\quad 2\le i<j\le p,
   \vphantom{\vrule height 11pt}\\
\nabla R_f(\partial_1^x,\partial_i^x,\partial_i^x,\partial_1^x;\partial_1^x)=\Theta_{;111}
   &\text{for}\quad 2\le i\le p\,.\vphantom{\vrule height 11pt}
\end{array}$$
Consequently after taking into account to normalize the basis for the tangent bundle suitably, we have
$$\alpha_f=\textstyle\frac{4(p-1)\Theta_{;111}^2}{(1+\Theta_{;11})^3}\,.$$
It is now clear that the metric $g_f$ will not be locally homogeneous for generic $\Theta$.

\hfill\qedbox

\begin{remark}\label{gronk3.3}\rm Let $\{\vec u_1,...,\vec u_p,\vec v_1,...,\vec v_p\}$ be a basis for a vector space
$V$ of dimension
$2p$. Define an innerproduct $(\cdot,\cdot)$ and an algebraic curvature tensor $R$ on $V$ whose non-zero entries are
$$(\vec u_i,\vec v_j)=\delta_{ij}\quad\text{and}\quad
R(\vec u_i,\vec u_j,\vec u_k,\vec u_l)=\delta_{il}\delta_{jk}-\delta_{ik}\delta_{jl}\,.$$ Then by Lemma \ref{gronk2.2},
$(V,(\cdot,\cdot),R)$ is a model for the metric and curvature tensor of all the manifolds
$(M,g_f)$ considered above. If we set $\Theta=0$, then
$$f_0=\textstyle\frac12\{x_1^2+...+x_p^2\}\,.$$
Since $\nabla R=0$,
$(M,g_{f_0})$ is a symmetric space and hence locally homogeneous. This shows that $(V,(\cdot,\cdot),R)$ is the
model for a symmetric space. Thus there exist pseudo-Riemannian manifolds which are not locally homogeneous whose
metric and curvature tensor is modeled on those of a symmetric space.\end{remark}


\begin{thebibliography}{AAA}

\bibitem{BK96}  E. Boeckx,, O. Kowalski, and L. Vanhecke,  {\bf Riemannian manifolds of
conullity two}, 
 World Scientific Publishing Co., Inc., River Edge, NJ, 1996.

\bibitem{BV98}  E. Boeckx and L. Vanhecke, {\it Curvature homogeneous unit tangent sphere bundles},
Publ. Math. {\bf 53}  (1998), 389-413.

\bibitem{BCR01} A. Bonome, P. Castro, and E. Garc\'{\i}a--R\'{\i}o, {\it Generalized Osserman four-dimensional
manifolds}, Classical Quantum Gravity {\bf 18} (2001), 4813--4822.

\bibitem{BCGHV98} A. Bonome, R. Castro, E. Garc\'{\i}a--R\'{\i}o, L. M. Hervella, R. V\'{a}zquez-Lorenzo,
{\it Nonsymmetric Osserman indefinite K\"{a}hler manifolds}, Proc. Amer. Math. Soc., {\bf 126}
 (1998), 2763--2769.

\bibitem{BM00} P. Bueken and M. Djori\'c, {\it Three-dimensional Lorentz metrics and
curvature homogeneity of order one}, Ann. Global Anal. Geom. {\bf 18} (2000), 85--103.

\bibitem{BV97} P. Bueken and L. Vanhecke, {\it Examples of curvature homogeneous Lorentz
metrics}, Classical Quantum Gravity {\bf 14} (1997), L93--L96.

\bibitem{CLPT90} M. Cahen, J. Leroy, M. Parker, F. Tricerri, and L. Vanhecke, {\it Lorentz
manifolds modelled on a Lorentz symmetric space}, J. Geom. Phys. {\bf 7} (1990), 571--581.

\bibitem{GKV02} E. Garc\'{\i}a--R\'{\i}o, D. Kupeli, and R. V\'azquez-Lorenzo, {\bf Osserman Manifolds in
Semi-Riemannian Geometry}, Lecture notes in Mathematics, Springer Verlag, (2002), ISBN 3-540-43144-6.

\bibitem{G02} P. Gilkey, {\bf Geometric Properties of Natural Operators Defined by the Riemann
Curvature Tensor}, World Scientific ISBN 981-02-4752-4 (2002).

\bibitem{GIZ02}  P. Gilkey, R. Ivanova, and T. Zhang, {\it Higher order Jordan
   Osserman pseudo-Riemannian manifolds}, Classical and Quantum
   Gravity, {\bf 19} (2002), 4543-4551; math.DG/0205269.

\bibitem{GIZ03} P. Gilkey, R. Ivanova, and T. Zhang, {\it Szab\'o Osserman IP Pseudo-Riemannian
manifolds}, Publ. Math. Debrecen 62 (2003), 387--401; math.DG/0205085.

\bibitem{IS95}  R. Ivanova and G. Stanilov,  {\it A skew-symmetric curvature
operator in Riemannian geometr}y, {\bf  Sympos. Gaussiana, Conf A}, ed. Behara, Fritsch,
and Lintz (1995), 391--395.

\bibitem{Ka00} I. Kath, {\it Killing spinors on pseudo-Riemannian manifolds},
Habilitation, Humboldt Universit\"at zu Berlin (2000).


\bibitem{K98} O. Kowalski, {\it On curvature homogeneous spaces}, 
Cordero, L. A. (ed.) et al., {\bf Proceedings of the workshop on recent topics in differential geometry},
Santiago de Compostela, Spain, July 16--19, 1997. Santiago de Compostela: Universidade de Santiago de
Compostela. Publ. Dep. Geom. Topolog\'a, Univ. Santiago Compostela. 89, 193-205 (1998).

\bibitem{KO99} O. Kowalski, B. Opozda, and Z. Vl\v sek, {\it Curvature homogeneity of affine connections on
two-dimensional manifolds}, Colloq. Math. {\bf 81} (1999), 123--139.

\bibitem{KTV91} O. Kowalski, Oldrich, F. Tricerri, and L. Vanhecke, {\it New examples of non-homogeneous
Riemannian manifolds whose curvature tensor is that of a Riemannian symmetric space},  C. R. Acad. Sci., Paris,
S\'er. I {\bf 311} (1990), 355-360.

\bibitem{KTV92} O. Kowalski, F. Tricerri and L. Vanhecke,
{\it Curvature homogeneous Riemannian manifolds},
 J. Math. Pures Appl., {\bf 71} (1992), 471--501.

\bibitem{O96} B. Opozda, {\it On curvature homogeneous and locally homogeneous affine connections}, Proc. Am.
Math. Soc. {\bf 124} (1996), 1889--1893.

\bibitem{Os97} R. Osserman, {\it Curvature in the eighties}, Amer. Math.
    Monthly, {\bf97}, (1990) 731--756.

\bibitem{Opr89} V. Oproiu, {\it Harmonic maps between tangent bundles}, Rend. Sem. Mat. Univ. Politec. Torino
{\bf 47} (1989),
47--55.

\bibitem{SV98} G. Stanilov and V. Videv, {\it Four dimensional
    pointwise Osserman manifolds}, Abh. Math. Sem. Univ.
    Hamburg, {\bf 68}, (1998),  1--6.

\bibitem{Sz91} Z. I. Szab\'o,  {\it A short topological proof for the symmetry of $2$ point homogeneous
spaces}, Invent. Math., {\bf 106}, (1991), 61--64.

\bibitem{T97} A. Tomassini, {\it Curvature homogeneous metrics on principal fibre bundles},
Ann. Mat. Pura Appl., IV. Ser. {\bf 172} (1997), 287--295.

\bibitem{T88} F. Tricerri, {\it Riemannian manifolds which have the same curvature as a homogeneous space, and a
conjecture of Gromov},  Riv. Mat. Univ. Parma, IV. Ser. {\bf 14} (1988), 91--104.

\bibitem{V91} L. Vanhecke, {\it Curvature homogeneity and related problems},
Recent topics in differential geometry, Proc. Workshop/Puerto de la Cruz/Spain 1990, Ser. Inf. 32 (1991),
103-122.

\end{thebibliography}
\end{document}